\documentclass[11pt]{amsart}

\makeatletter
\usepackage{amssymb}
\usepackage{latexsym}
\usepackage{amsbsy}
\usepackage{amsfonts}
\usepackage[notref,notcite]{showkeys}

\def\marginpar#1{\ignorespaces}

\newtheorem{theorem}{Theorem}
\newtheorem{proposition}[theorem]{Proposition}
\newtheorem{lemma}[theorem]{Lemma}
\newtheorem{corollary}[theorem]{Corollary}

\theoremstyle{definition}
\newtheorem{remark}[theorem]{Remark}

\numberwithin{equation}{section}
\numberwithin{theorem}{section}

\def\AArm{\fam0 \rm}%
\newdimen\AAdi%
\newbox\AAbo%
\def\AAk#1#2{\setbox\AAbo=\hbox{#2}\AAdi=\wd\AAbo\kern#1\AAdi{}}%

\newcommand{\BBone}{{\ensuremath{{\AArm 1\AAk{-.8}{I}I}}}}

\def\eqref#1{(\ref{#1})}
\def\eqlabel#1{\def\@currentlabel{#1}}

\def\formula#1{\def\@tempa{#1}\let\@tempb\theequation\def\theequation{%
\hbox{#1}}\def\@currentlabel{(\theequation)}$$}
\def\endformula{\leqno\hbox{(\@tempa)}$$\@ignoretrue\let\theequation\@tempb}

\def\given{\hskip5\p@\relax\vrule\@width.4\p@\hskip5\p@\relax}

\newcommand{\open}[1]{%
\par\normalfont\topsep6\p@\@plus6\p@\trivlist\item[\hskip\labelsep\itshape#1%
\@addpunct{.}]\ignorespaces}

\DeclareRobustCommand{\close}[1]{%
  \ifmmode 
  \else \leavevmode\unskip\penalty9999 \hbox{}\nobreak\hfill
  \fi
  \quad\hbox{$#1$}}

\newlength{\toskip}

\def\LL{\mathcal L}

\def\<{\langle}
\def\>{\rangle}

\def \R {{\mathbb R}}

\def \P {{\mathbb P}}

\def \L {{\mathbb L}}

\def \Var {\textrm{Var}}
\def \Ent {\textrm{Ent}}

\makeatother

\begin{document}
\date{\today}

\title[Weighted logarithmic Sobolev inequality]{Some remarks on weighted logarithmic Sobolev inequality}

 \author[P. Cattiaux]{\textbf{\quad {Patrick} Cattiaux $^{\spadesuit}$ \, \, }}
\address{{\bf {Patrick} CATTIAUX},\\ Institut de Math\'ematiques de Toulouse. CNRS UMR 5219. \\
Universit\'e Paul Sabatier,
\\ 118 route
de Narbonne, F-31062 Toulouse cedex 09.} \email{cattiaux@math.univ-toulouse.fr}

 \author[A. Guillin]{\textbf{\quad {Arnaud} Guillin $^{\diamondsuit}$}}
\address{{\bf {Arnaud} GUILLIN},\\ Laboratoire de Math\'ematiques, CNRS UMR 6620, Universit\'e Blaise Pascal, avenue des Landais 63177 Aubi\`ere.}
\email{guillin@math.univ-bpclermont.fr}

 \author[L-M. Wu]{\textbf{\quad {Li-Ming} Wu $^{\diamondsuit}$}}
\address{{\bf {Li-Ming} WU},\\ Laboratoire de Math\'ematiques, CNRS UMR 6620, Universit\'e Blaise Pascal, avenue des Landais 63177 Aubi\`ere.}
\email{wuliming@math.univ-bpclermont.fr}

\maketitle
 \begin{center}

 \textsc{$^{\spadesuit}$  Universit\'e de Toulouse}
\smallskip

\textsc{$^{\diamondsuit}$ Universit\'e Blaise Pascal}
\smallskip

 \end{center}

\begin{abstract}
We give here a simple proof of weighted logarithmic Sobolev inequality, for example for Cauchy
type measures, with optimal weight, sharpening results of Bobkov-Ledoux \cite{BLweight}. Some
consequences are also discussed.
\end{abstract}
\bigskip

\bigskip

\bigskip

\section{Introduction}

In a recent paper, Bobkov and Ledoux \cite[Th. 3.4]{BLweight} proved that for a generalized Cauchy
measure on $\R^n$, i.e.
$$d\nu_\beta(x)=\frac1Z(1+|x|^2)^{-\beta}dx$$
for $\beta>n/2$, the following weighted logarithmic Sobolev
inequality holds, provided $\beta \geq (n+1)/2$: for any smooth
bounded $f$
$$\Ent_{\nu_\beta}(f^2)=\nu_\beta\left(f^2\log\left( \frac{f^2}{\nu_\beta(f^2)}\right)\right)
\le \frac{1}{\beta-1}\int|\nabla f(x)|^2(1+|x|^2)^2d\nu_\beta(x).$$ Simple test functions however
indicate that the weight $(1+|x|^2)^2$ is not optimal: one hopes $(1+|x|^2)\log(e+|x|^2)$ and that
is what we will recover (with somewhat less precise constants).

\medskip

It will be thus our purpose to prove inequalities of the type
$$\Ent_\mu(f^2)\le c\int |\nabla f|^2\omega d\mu$$
for some weight $\omega\ge1$, and more generally weighted $F$-Sobolev inequalities with more
general $F$'s replacing the logarithm.
\smallskip

The (in a particular sense) case of weighted Poincar\'e inequalities
is studied in \cite{BLweight} for Cauchy type measures and in
\cite{CGGR} in more general situations. Consequences in terms of
concentration of measure or isoperimetry are described in details in
the latter reference.

It should also be interesting to look at weights that go to 0 at infinity (instead of weights
bounded by 1 from below). Part of the results in \cite{CGGR} and in the present paper extend to
this situation.
\medskip

Our strategy will be the following:
\begin{enumerate}
\item consider a dual form of the weighted logarithmic Sobolev inequality (or more generally
$F$-Sobolev inequality): the Super weighted Poincar\'e inequality; \item use Lyapunov condition to
prove these Super weighted Poincar\'e inequalities; \item show that these Super weighted
Poincar\'e inequalities are equivalent to  weighted F-Sobolev inequality (and in particular
weighted logarithmic Sobolev inequality).
\end{enumerate}

Let us then introduce the so called Super weighted Poincar\'e inequality for a probability measure
$\mu$, in a simple context, namely when the underlying Carr\'e du champ is in fact the square
length of the gradient. It is inspired from the pioneering works on Super Poincar\'e inequality by
Wang \cite{w00}. Given a weight $\omega$ larger than 1, we say that $\mu$ satisfies a Super
weighted Poincar\'e inequality if for all $f$ smooth and bounded, there exists  a non-increasing
function $\beta_\omega$ such that for all $s>0$
\begin{equation}\label{SwPI}
\int f^2d\mu\le s \int |\nabla f|^2\omega d\mu+\beta_\omega(s)(\mu(|f|))^2 \, .
\end{equation}
When $\omega=1$, it is the usual Super Poincar\'e inequality which describes properties of the
measure stronger than the usual Poincar\'e inequality. If we add some additional weight $\omega$
(tending to infinity as $|x|\to\infty$ for example) we will be able to give an inequality adapted
to measures ``above'' and ``below'' Poincar\'e, being even able to play between the weight and
$\beta$.

Weighted Poincar\'e inequalities have been recently investigated by
Bobkov-Ledoux \cite{BLweight} in particular for their interest in
deviation inequalities, and by Cattiaux and al \cite{CGGR} showing
their link with weak Poincar\'e inequalities and isoperimetric
inequalities. They have been also intensively studied, in a converse
form,  in PDE theory to establish exponential convergence to
equilibrium for fast diffusion equations (see \cite{DMC,blanchet}).
In parallel, Cattiaux and al \cite{CGWW} have studied Super
Poincar\'e inequalities using Lyapunov conditions (see also
\cite{BBCG,BCG}). We will combine here these two approaches to study
these Super weighted Poincar\'e inequalities.

\section{Results and examples}

\subsection{A Lyapunov condition for Super weighted Poincar\'e inequality}

Lyapunov conditions appeared a long time ago in relation with the
problem of convergence to equilibrium for Markov processes, see
\cite{MT,MT2,MT3,DFG} and references therein. They also have been
used to study large and moderate deviations for empirical
functionals of Markov processes (see Donsker-Varadhan
\cite{DV3,DV4}, Kontoyaniis-Meyn \cite{KM1,KM2}, Wu \cite{wu1},
Guillin \cite{G1,G2},...) Their use to provide functional
inequalities has been very recently deeply investigated  with some
success: Lyapunov-Poincar\'e inequalities \cite{BCG}, Poincar\'e
inequalities \cite{BBCG}, transportation inequalities \cite{CGW},
Super Poincar\'e inequalities \cite{CGWW}, weighted and weak
Poincar\'e inequalities \cite{CGGR} (also see the recent survey
\cite{CGgre}). We will take advantage of the approach of these last
two papers to build our main results, but let us first describe our
framework.
\medskip

Let $E$ be some Polish state space, $\mu$ a probability measure and
a $\mu$-symmetric diffusion semigroup $P_t$ with generator $L$ on
$L^2(E,\mu)$. The main assumption on $L$ is that there exists some
algebra $\mathcal{A}$ of bounded and uniformly continuous functions,
containing constant functions, which is  in the domain of $L$ in the
graph norm of $\LL$ on $L^2(\mu)$. It enables us to define a
``carr\'e du champ'' $\Gamma$, {\it i.e.} for $f, g \in
\mathcal{A}$, $L(fg)=f Lg + g Lf + 2 \Gamma(f,g)$. We will also
assume that $\Gamma$ is a derivation (in each component), {\it i.e.}
$\Gamma(fg,h)=f\Gamma(g,h) + g \Gamma(f,h)$, i.e. we are in the
standard ``diffusion'' case in \cite{bakry} and we refer to the
introduction of \cite{cat4} for more details. For simplicity we set
$\Gamma(f)=\Gamma(f,f)$. Also, since $L$ generates a diffusion, we
have the following chain rule formula $\Gamma(\Psi(f),
\Phi(g))=\Psi'(f) \Phi'(g) \Gamma(f,g)$.

In particular if $E=\R^n$, $\mu(dx)= p(x) dx$ and $L=\Delta + \nabla
\log p\cdot\nabla$, we may consider the algebra generated by
$C^\infty$ functions with compact support and the constant
functions, as the interesting subalgebra $\mathcal{A}$, and then
$\Gamma(f,g)=\nabla f \cdot \nabla g$.
\medskip

Now we define the notion of $\phi$-Lyapunov function. Let $W\geq 1$ be a smooth enough function on
$E$ and $\phi$ be a $\mathcal{C}^1$ positive increasing function defined on $\R^+$. We say that
$W$ is a $\phi$-Lyapunov function if there is a family of increasing sets $(A_r)_{r\ge0} \subset
E$ such that $\bigcup_r A_r = E$ (we say that the family $A_r$ is exhausting) and some $b \geq 0$
such that for some $r_0>0$
\begin{equation}
LW \, \le \,  -\phi(W) \, + \, b \, \BBone_{A_{r_0}} \, .\label{lyap}
\end{equation}
This latter condition is sometimes called a ``drift condition'' but we prefer to call it Lyapunov
condition.  One has very different behavior depending on $\phi$: if $\phi$ is linear then a
Poincar\'e inequality is valid, whereas when $\phi$ is super-linear (or more generally in the form
$\phi\times W$ where $\phi$ tends to infinity ) we have stronger inequalities (Super Poincar\'e,
ultracontractivity...), and finally if $\phi$ is sub-linear  we are in the regime of weak
Poincar\'e inequalities. We will cover setting in both weak and super Poincar\'e inequalities
playing with the weight function.

\medskip

We are now in position to state our main theorem:

\begin{theorem}\label{mainTH}
Assume that $L$ satisfies a Lyapunov condition (\ref{lyap}), that $\mu$ satisfies some local Super
 Poincar\'e inequality, i.e. there exists $\beta_{loc}$ decreasing in $s$ (for all $r$) such that
$\forall s>0$
\begin{equation}\label{SPloc}
\int_{A_r} f^2d\mu\le s\int \Gamma(f)d\mu+\beta_{loc}(r,s)\left(\int_{A_r}|f|d\mu\right)^2 \, .
\end{equation}
We also introduce some $\psi:[1,\infty[\to[1,\infty[$ which is
increasing and such that $$0 < (\phi/\psi)'(W)\le1 \, .$$ We finally
assume that $G(r):=1/(\inf_{A_r^c}\psi(W))$ goes to 0 as $r$ goes to
infinity.

Then $\mu$ satisfies a Super weighted Poincar\'e inequality, i.e. for all $s>0$
\begin{equation}
\int f^2d\mu\le 2s \int \, \frac{\Gamma(f)}{\left(\frac{\phi}{\psi}\right)'(W)} \, \,
d\mu+\tilde\beta(s)\left(\int |f|d\mu\right)^2 \label{result}
\end{equation}
where
$$\tilde\beta(s)=c_{r_0}\,\beta_{loc}(G^{-1}(s),s/c_{r_0})$$
$G^{-1}(s)=\inf\{t>0; G(t)>s\}$ is the right inverse of $G$ and
$$c_{r_0}=1+b\frac{\sup_{A_{r_0}}(\psi/\phi)(W)}{\inf_{A_{r_0}^c}\psi(W)}.$$
\end{theorem}

\begin{remark}
In fact it is of course sufficient to verify some local Super weighted Poincar\'e inequality,  but
as the weight is usually bounded on the subset $A_r$ considered, they are equivalent (up to the
constants involved). And even, playing with $r$, as the weight is supposed to be greater than 1
they are implied by the local Super Poincar\'e inequalities as used here.
\end{remark}

\begin{remark}
In the particular case where $\Gamma(f,g)=\nabla f\cdot\nabla g$ one
can take more general Lyapunov condition, namely $\phi(W)$ may be
replaced by $\phi\times W$ for some functional $\phi$ and the same
for $\psi$ appearing in the theorem. The modifications are
straightforward but give hard to read result, and we let then the
details for people needing such a framework.
\end{remark}

\begin{remark}
In practice, $A_r$ will often be level sets of the Lyapunov function
$W$  or balls of radius $r$. The local Super Poincar\'e inequality
will then be obtained by perturbation of the Super weighted
Poincar\'e inequality on balls for the underlying (Lebesgue)
measure.
\end{remark}

\begin{proof}
Let us begin with quite easy estimates: for $r\ge r_0$
\begin{eqnarray*}
\int f^2d\mu&=&\int_{A_r} f^2d\mu+\int_{A_r^c} f^2d\mu\\
&=&\int_{A_r} f^2d\mu+\int_{A_r^c} \frac{\psi(W) \, \phi(W)}{\psi(W) \, \phi(W)}f^2d\mu\\
&\le&\int_{A_r} f^2d\mu+\frac1{\inf_{A_r^c}\psi(W)}\int f^2 \frac{\psi(W)}{\phi(W)} \, \phi(W) \,  d\mu\\
&\le& \int_{A_r} f^2d\mu+b \, \frac{\sup_{A_{r_0}}\left(\frac{\psi}{\phi} \, (W)\right)}{\inf_{A_{r}^c}\psi(W)} \, \int_{A_{r_0}} f^2d\mu\\
&+& \frac1{\inf_{A_r^c}\psi(W)}\int \frac {-LW}{\frac{\phi}{\psi}(W)} \, f^2 \, d\mu\\&\le&
\left(1 + b \, \frac{\sup_{A_{r_0}}\left(\frac{\psi}{\phi} \,
(W)\right)}{\inf_{A_{r}^c}\psi(W)}\right) \, \int_{A_r} f^2d\mu + \frac1{\inf_{A_r^c}\psi(W)}\int
\frac {-LW}{\frac{\phi}{\psi}(W)} \, f^2 \, d\mu \, .
\end{eqnarray*}
Applying Lemma \ref{lem:philyapounov} below to the second term, the
local Super Poincar\'e inequality and the fact that
$(\phi/\psi)'(W)\le1$ to the first, we get
\begin{eqnarray*}
\int f^2d\mu&\le&\left(s\left(1+b\frac{\sup_{A_{r_0}}(\psi/\phi)(W)}{\inf_{A_r^c}\psi(W)}\right)+
\frac1{\inf_{A_r^c}\psi(W)}\right)\int \, \frac{\Gamma(f)}{\left(\frac{\phi}{\psi}\right)'(W)} \, d\mu\\
&&+\beta_{loc}(r,s) \left(1+b\frac{\sup_{A_{r_0}}(\psi/\phi)(W)}{\inf_{A_r^c}\psi(W)}\right)
\left(\int|f|d\mu\right)^2 \, .
\end{eqnarray*}
Recall now
$$c_{r_0}=1+b\frac{\sup_{A_{r_0}}(\psi/\phi)(W)}{\inf_{A_{r_0}^c}\psi(W)}$$
and $\tilde s=sc_{r_0}$ so that, since $A_r^c$ is decreasing in $r$, the last inequality furnishes
$$\int f^2d\mu\le \left(\tilde s+G(r)\right) \int \, \frac{\Gamma(f)}{\left(\frac{\phi}{\psi}\right)'(W)} \, d\mu +
 \beta_{loc}(r,\tilde s/c_{r_0})c_{r_0} \left(\int|f|d\mu\right)^2.$$
Choose now $r=G^{-1}(\tilde s)$ to conclude.
\end{proof}

One crucial element of the proof above was the following lemma borrowed from \cite{CGGR}  whose
proof is reproduced here for completeness (showing also the necessity for $L$ to be a diffusion)
\begin{lemma} \label{lem:philyapounov}
 Let $\psi : \mathbb{R}^+ \to \mathbb{R}^+$ be a $\mathcal{C}^1$ increasing function.
 Then, for any $f\in \mathcal{A}$ and any positive $h\in D(\mathcal E)$,
$$
\int  \, \frac{-Lh}{\psi(h)} \, f^2 \, d\mu  \leq  \int \frac { \, \Gamma(f)}{ \psi'(h)} d\mu
$$
\end{lemma}

\begin{proof}
Since $L$ is $\mu$-symmetric, using that $\Gamma$ is a derivation and the chain rule formula, we
have
\begin{eqnarray*}
\int  \, \frac{-Lh}{\psi(h)} \, f^2 \, d\mu
& = & \int \, \Gamma \left( h, \frac{f^2}{\psi(h)} \right) \, d\mu
=
\int \left(\frac{2 \, f \,
\Gamma(f,h)}{\psi(h)} \, - \, \frac{f^2 \psi'(h) \Gamma(h)}{\psi^2(h)}\right) d\mu \, .
\end{eqnarray*}
Since $\psi$ is increasing and according to Cauchy-Schwarz inequality we get
\begin{eqnarray*}
\frac{f \, \Gamma(f,h)}{\psi(h)}
& \leq &
\frac{f \sqrt{\Gamma(f) \Gamma(h)}}{\psi(h)}
 =
\frac{ \sqrt{\Gamma(f)} }{ \sqrt{\psi'(h)} } \cdot
\frac{f \sqrt{\psi'(h) \Gamma(h)}}{\psi(h)}
\\
&\leq &
\frac{1}{2} \frac{\Gamma(f)}{\psi'(h)} + \frac{1}{2} \, \frac{f^2 \psi'(h)\,
\Gamma(h)}{\psi^2(h)}.
\end{eqnarray*}
The result follows.
\end{proof}

\subsection{Equivalence with weighted $F$-Sobolev inequality}

Let $F$ be a continuous function, such that
$\sup_{0<r<1}|rF(r)|<\infty$, $F(1)=0$ and
$\lim_{x\to+\infty}F(x)=+\infty$. We will say that the probability
measure $\mu$ satisfies a defective weighted $F$-Sobolev inequality,
with constants $C_1$ and $C_2$, and weight $\omega$, if for all
smooth and bounded $f$ with $\mu(f^2)=1$
$$\int f^2F(f^2)d\mu\le C_1\int \Gamma(f) \, \omega d\mu+C_2.$$ Notice that, modifying if necessary
the constant $C_2$ we may replace $F$ by $F_+$. This inequality will be called tight, or simply a
weighted $F$-Sobolev inequality if $C_2=0$.
\medskip

When $\omega=1$, it is known that if $\mu$ satisfies a defective $F$-Sobolev inequality and a
Poincar\'e inequality, and with some (slight) additional assumptions on $F$, then $\mu$ satisfies
a (tight) $F$-Sobolev inequality. The case $F=\log$ is known as Rothaus lemma, and the previous
general result is obtained in \cite{BCR1} lemma 9 and Theorem 10.

The reader will easily check that the proofs in \cite{BCR1} extend to the weighted case, i.e. a
weighted Poincar\'e inequality (with weight $\omega$) and a weighted defective $F$-Sobolev
inequality (with the same $\omega$) imply a tight weighted $F$-Sobolev inequality, under the same
assumptions than in \cite{BCR1} lemma 9. These assumptions are satisfied when $F(x)=\log_+(x)$
(see remark 15 in \cite{BCR1}). We thus have that a weighted $\log$-Sobolev inequality implies a
weighted $\log_+$-Sobolev inequality, and together with a weighted Poincar\'e inequality implies a
tight weighted $\log_+$-Sobolev inequality, hence a tight weighted $\log$-Sobolev inequality.

We shall use this line of reasoning in various situations below, without mentioning it explicitly.
\medskip

Now let us make a simple remark: if in the Super weighted Poincar\'e inequality, we assume
moreover that $\beta_\omega$ tends to a constant smaller than 1 as $s\to\infty$ (which is quite a
very weak hypothesis), the Super weighted Poincar\'e inequality implies a weighted Poincar\'e
inequality. Indeed applying \eqref{SwPI} with $f= g -\mu(g)$ we get
$$(1-\beta_\omega(s)) \, \Var_\mu(g) \leq s \, \int \, \Gamma(g) \, \omega  \, d\mu \, ,$$ thanks to Cauchy-Schwarz inequality,
and the
result follows taking a large enough $s$ for the left hand side to be positive.
\medskip

The next proposition is adapted from the works of Wang \cite{w00} and Theorems 3.3.1 and 3.3.3 in
\cite{Wbook}. We include its proof for the sake of completeness.

\begin{proposition}\label{eqFsob}

\begin{enumerate}
\item If $\mu$ satisfies a defective weighted $F$-Sobolev inequality with constants $C_1$, $C_2$,
then there exist $c_1,c_2$ such that for all smooth bounded functions $f$ and $\forall s>0$
$$\int f^2 d\mu \le s\int \Gamma(f)\omega d\mu +c_1 F^{-1}(c_2(1+1/s))\mu(|f|)^2$$
where $F^{-1}(s)=\inf\{r\ge0;~F(r)\ge s\}$. \item If $\mu$ satisfies a Super weighted Poincar\'e
inequality
$$\int f^2 d\mu \le s \int \Gamma(f) \omega d\mu +\beta_\omega(s)\mu(|f|)^2$$
then $\mu$ satisfies a defective weighted $F$-Sobolev inequality with
$$F(r)=\frac{c_1(\epsilon)}{r}\int_0^r\xi(\epsilon t)dt-c_2(\epsilon)$$
for all $0<\epsilon<1$, where $c_1(\epsilon)$ and $c_2(\epsilon)$ are some constants, and
$$\xi(t)=\sup_{r>0}\left(\frac 1r -\frac{\beta_\omega(r)}{rt}\right),$$
where $\beta_\omega^{-1}(t)=\inf\{r\ge0;\beta_\omega(r)\le t\}.$
\end{enumerate}
\end{proposition}

\begin{proof}
(1). As said before we may assume that $F\ge 0$, enlarging $C_2$ if necessary. Pick $f$ with
$\mu(|f|)=1$. For all $r,t,a>0$, it holds
$$rt\le rF(r^2/a)+t\sqrt{aF^{-1}(t)} \, .$$
We choose $a=\mu(f^2)$, $r=|f|$ and multiply the previous inequality by $|f|$, i.e. $$t \, f^2
\leq f^2 \, F(f^2/\mu(f^2)) + |f| t \sqrt{\mu(f^2) \, F^{-1}(t)} \, .$$ Integrating this
inequality with respect to $\mu$ yields
$$\mu(f^2F(f^2/\pi(f^2)))\ge t\mu(f^2)-t\sqrt{\mu(f^2)F^{-1}(t)}$$
and using the defective weighted $F$-Sobolev inequality :
$$(t-C_2)\mu(f^2)-t\sqrt{\mu(f^2)F^{-1}(t)} - C_1 \int \Gamma(f) \, \omega d\mu \le 0 \, .$$
Hence, for $t>C_2$,
$$\mu(f^2)\le \frac{2C_1}{t-C_2}\int \Gamma(f) \omega d\mu + \frac{t^2F^{-1}(t)}{(t-C_2)^2} \, ,$$
and we write $r=2C_1/(t-C_2)$ to conclude.\\
(2).  The second part of the proof is inspired by capacity/measure criteria.\\
Pick $f$ with $\mu(f^2)=1$ and $\delta>1$, consider $A_n=\{\delta^{n+1}>f^2\ge\delta^n\}$ and
$$f_n=(|f|-\delta^{n/2})_+\wedge(\delta^{(n+1)/2}-\delta^{n/2}).$$ Apply now the Super weighted Poincar\'e inequality to $f_n$,
$$\mu(f_n^2)\le r\mu(\Gamma(f)\omega\BBone_{A_n})+\beta(r)\mu(f_n)^2\le  r\mu(\Gamma(f)\omega\BBone_{A_n})+
\beta_\omega(r)\mu(f^2\ge \delta^n)\mu(f_n^2)$$ and since $\mu(f^2\ge \delta^n)\le 1/\delta^n$, we
get
\begin{eqnarray*}
\mu(\Gamma(f)\omega)&\ge& \sum_{n\ge 0}\mu(\Gamma(f)\omega\BBone_{A_n})\\
&\ge& \sum_{n\ge0}\xi(\delta^n)\mu(f^2_n)\\
&\ge&\sum_{n\ge 0} \xi(\delta^n)\mu(f^2\ge\delta^{n+1})(\delta^{(n+1)/2}-\delta^{n/2})^2\\
&\ge&\frac{(\sqrt{\delta}-1)^2}{1-\delta^{-1}}\sum_{n\ge0}\int_{\delta^{n-1}}^{\delta^n}\xi(t)\mu(f^2\ge\delta^2t)dt\\
&\ge&c_1\int_0^\infty \xi(t))\mu(f^2\ge\delta^2t)dt-c_2\\
&\ge&c_3\pi(f^2F(f^2))-c_2
\end{eqnarray*}
which is what is needed.
\end{proof}

Using this result one sees that if a Super weighted (with weight $\omega$) Poincar\'e  inequality
is valid with $\beta_\omega(s)= s^{-N} \, e^{c(1+1/s)}$ then a ($\omega$) weighted  logarithmic
Sobolev inequality is valid. In the preceding subsection we have presented conditions to verify
Super weighted Poincar\'e inequalities, we only have now to validate them through examples. It
will be the purpose of the next subsection.

\subsection{Examples}

We consider here the $\R^n$ situation with $d \mu(x)=p(x) dx$ and
$L=\Delta + \nabla \log p\cdot\nabla$, where $p$ is smooth enough
and positive, and $\cdot$ is the euclidean inner product. Recall the
following elementary lemma whose proof can be found in \cite{BBCG}.
This lemma will be helpful to deal with $\kappa$-concave measures.

\begin{lemma}\label{lemfranck}
If $V$ is convex and $\int e^{-V(x)} \, dx < +\infty$, then
\begin{itemize} \item[(1)] \quad for all $x$, $x\cdot\nabla V(x) \geq V(x) - V(0)$, \item[(2)] \quad
there exist $\delta >0$ and $R>0$ such that for $|x|\geq R$,
$V(x)-V(0) \geq \delta \, |x|$.
\end{itemize}
\end{lemma}

Another helpful result is the following result concerning the validity of a Super Poincar\'e
inequality for Lebesgue measures on balls: for all $r>0$ denote by $B(0,r)$ the euclidean ball in
$\R^n$. Then there exists $c_n$ such that for all smooth $f$ and all $s>0$,
\begin{equation}
\label{SPIballs}
\int_{B(0,r)}f^2dx\le s \int_{B(0,r)}|\nabla f|^2dx+c_n(1+s^{-n/2})\left(\int_{B(0,r)}|f|dx\right)^2.
\end{equation}
Such an inequality will be particularly efficient when dealing with radial type measures, as
perturbation argument to get the local Super Poincar\'e inequality will be easy to do.

Indeed we immediately obtain
\begin{eqnarray}
\label{muSPIballs} \int_{B(0,r)}f^2d\mu &\le& s \int_{B(0,r)}|\nabla f|^2d\mu
\\ &+& c_n\left(1+\left(\frac{ s \, \inf_{B(0,r)} p}{\sup_{B(0,r)}p}\right)^{-n/2}\right) \,
\left(\frac{\sup_{B(0,r)} p}{\inf^2_{B(0,r)}p}\right) \, \left(\int_{B(0,r)}|f|d\mu\right)^2 \, .
\nonumber
\end{eqnarray}

For more general type of measures, it is not so difficult to get local inequalities for level sets
of the potential, see \cite[Prop. 3.6]{CGWW}.

\subsubsection{Cauchy type measures}

Let $d\mu(x) = (V(x))^{-(n+\alpha)} \, dx$ for some positive convex function $V$ and some $\alpha>
0$. Let us begin by establishing a Lyapunov condition:

\begin{lemma}\label{lem:convex}
Let $L= \Delta - (n+\alpha) (\nabla V/V) \nabla$ with $V$ convex and
$\alpha>0$. Then, there exists $k\in (2,\alpha+2)$,  $b, R>0$ and
function $W \geq 1$ such that
$$
LW \leq - \phi(W) + b \BBone_{B(0,R)}
$$
with $\phi(u)=c u^{(k-2)/k}$ for some constant $c>0$.  Furthermore, one can choose
$W(x)=|x|^k$ for $x$ large.
\end{lemma}

\begin{proof}
Let $L= \Delta - (n+\alpha) (\nabla V/V) \nabla$ and choose $W \geq 1$ smooth, satisfying
$W(x)=|x|^k$ for $|x|$ large enough and $k>2$ that will be chosen later. For $|x|$ large enough we
have
$$
LW(x) = k \, (W(x))^{\frac{k-2}{k}} \, \left(n+k-2 - \frac{(n+\alpha) \, x.\nabla
V(x)}{V(x)}\right) \, .
$$
Using (1) in Lemma \ref{lemfranck} (since $V^{-(n+\alpha)}$ is integrable, $e^{-V}$ is also
integrable)  we have
$$
n+k-2 - \frac{(n+\alpha) \,
x.\nabla V(x)}{V(x)} \leq k-2 - \alpha + (n+\alpha) \frac{V(0)}{V(x)} \, .
$$
Using (2) in Lemma
\ref{lemfranck} we see that we can
choose $|x|$ large enough for $\frac{V(0)}{V(x)}$ to be less than $\varepsilon$, say
$|x|>R_\varepsilon$. It remains to choose $k>2$ and $\varepsilon>0$ such that
$$
k + n \varepsilon -2 - \alpha (1-\varepsilon)\leq -\gamma
$$
for some $\gamma>0$. We have shown that, for $|x|>R_\varepsilon$,
$$
LW \leq - k \gamma \phi(W),
$$
with $\phi(u) =u^{\frac{k-2}{k}}$ (which is increasing since $k>2$). A compactness argument achieves the proof.
\end{proof}

Consider now the case studied in \cite{BLweight} of the (generalized) Cauchy measure:
$$p(x)=Z_\beta^{-1}(1+|x|^2)^{-\beta},\qquad \beta>n/2.$$

Lemma \ref{lem:convex} gives us a Lyapunov conditions. Using \eqref{muSPIballs} we get  local
Super Poincar\'e inequalities
\begin{eqnarray*}
\int_{B(0,R)}f^2d\mu &\le& s \int_{B(0,R)}|\nabla f|^2d\mu
\\ &+& c_n\left(1+s^{-n/2} \, (1+R^2)^{\beta \, n/2}\right) (1+R^2)^{2\beta} \, Z_\beta \left(\int_{B(0,R)}|f|d\mu\right)^2.
\end{eqnarray*}

Choose now  $\psi(v)=\log(v)$ for large $v$ (and $\psi$ smooth), Theorem \ref{mainTH} together
with Proposition \ref{eqFsob} thus furnishes (up to local modifications i.e for large $|x|$'s for
example)
$$\phi(u)=u^{k-2/k} \, , \, \psi(u) = \log (u) \, , \, W(x)=|x|^k \, , \, (\psi(W))(x)=k \log
|x|$$ hence $$G(r)= \frac{1}{k \, \log r} \, , \, G^{-1}(s) = e^{1/k s}$$ so that
$$\left(\frac{\phi}{\psi}\right)'(u) \sim \frac{c}{u^{2/k} \, \log u} \, , \, \omega(x) \sim
\left(\frac{1}{(\phi/\psi)'(W)}\right)(x) \sim c \, |x|^2 \, \log |x|$$ and finally for small $s$
$$\beta_\omega(s) \sim s^{- \, n/2} \, e^{c/s} \, .$$ We have thus obtained

\begin{corollary}\label{corcauchy}
Cauchy measures $\mu(dx)=Z_\beta^{-1}(1+|x|^2)^{-\beta}$ for $\beta>n/2$ verify the  following
weighted logarithmic Sobolev inequality: there exists $C=C(\beta,n)$ such that for all smooth
bounded function $f$
$$\Ent_\mu(f^2)\le C\int |\nabla f(x)|^2\,(1+|x|^2)\log(e+|x|^2)d\mu(x) \, .$$
\end{corollary}

We then obtain the correct order of magnitude of the weight in this  inequality, compared to
\cite[Th.3.4]{BLweight}. However it has to be noted that we are loosing the pretty expression of
the constant in front of the weighted energy. Note that in dimension 1, Barthe-Zhang \cite{BZ}
obtained the same weight.

\subsubsection{Exponential measure}

We will look at the  exponential measure
$$\nu(dx)= Z_n^{-1} \, e^{-|x|}dx.$$
It is well known that the exponential measure satisfies a Poincar\'e inequality.  It is also easy
to see that considering $W(x)=e^{a|x|}$ for $|x|\ge R$, we get if $a<1$ for $R$ large enough
$$LW(x)=a\left(\frac{n-1}{|x|} +a-1\right)W(x)\le -\lambda W+b\BBone_{B(0,R)}$$
and thus the Lyapunov condition.

Using \eqref{muSPIballs} with the choice $\psi(v)=\log(v)$ for large $v$ (and $\psi$ smooth), we
get
\begin{corollary}
The exponential measure $\nu$ satisfied the following weighted logarithmic Sobolev inequality:
there exists $C=C(\beta,n)$ such that for all smooth bounded function $f$,
$$\Ent_\mu(f^2)\le C\int |\nabla f(x)|^2 \, (1+|x|) \, d\mu(x).$$
\end{corollary}
As a comparison, let us recall a result of Bobkov-Ledoux \cite[Eq. (1.6)]{BL97} which  states that
for the one sided exponential $\tilde\nu$ (in dimension one)
$$\Ent_{\tilde\nu}(f^2)\le 4\int x(f'(x))^2d\tilde\nu.$$
We then recover in any dimension their result directly
(they can only use tensorization to get $n$-dimensional version of this inequality) and may extend it to other potential.
\medskip

\begin{remark}
Actually the proof above covers a very large class of measures satisfying a Poincar\'e inequality,
namely measures $\mu(dx) = e^{-V} dx$ such that $V \to +\infty$ as $|x| \to +\infty$ and
satisfying the following condition $$\textrm{ there exists $0<a<1$ such that } \liminf_{|x|\to
+\infty}\left(a |\nabla V|^2 - \Delta V\right) = B > 0 \, .$$ Indeed in this case we have
$\phi(u)=\lambda u$ (for some $\lambda >0$) and $W=e^{AV}$ for some well chosen positive constant
$A$.

Choosing again $\psi(u)=\log u$ for large $u$'s we obtain the weight $\omega(x)=|x|$ for large
$|x|$'s. If we assume in addition that there exists some constant $c>0$ such that for all $R$ and
all $x$ such that $|x|=R$, $$c \, \sup_{|y|=R} V(y) \le V(x) \le \frac 1c \, \inf_{|y|\ge R} V(y)
\, ,$$ it is not difficult to see that $G^{-1}(s) \sim (\bar V)^{-1}(1/s)$ where $\bar V(R) =
\inf_{|y|\ge R} V(y)$ is increasing. Using \eqref{muSPIballs} again we obtain that
$\beta_\omega(s) \sim \exp(C/s)$ hence the same weighted logarithmic Sobolev inequality as in the
previous corollary.

We do not know whether this is true for any measure satisfying the Poincar\'e inequality. Indeed
we know that there exists some Lyapunov function $W$ yielding a linear $\phi$, but we do not know
in full generality how to compare $W$ and the potential $V$, so that we cannot give an explicit
formula for $\beta_\omega$. \hfill $\diamondsuit$
\end{remark}

\section{Properties and Applications}

\subsection{Concentration of measure}

We will present here two different approaches to get concentration inequalities. The first one,
due to Bobkov-Ledoux \cite{BLweight} uses controls on the growth of moments. As we obtain optimal
weight by our approach, we will compare on some examples what are the implications of these better
controls. The other one is based on the derivation of a suitable transportation cost information
inequality following the approach of \cite{BGL} based on Hamilton-Jacobi equation.

\subsubsection{Growth of moments and Deviation inequality.}
We briefly recall here the main results concerning concentration inequality obtained by
Bobkov-Ledoux \cite[Th. 4.1,Cor. 4.2]{BLweight} and present their main result
\begin{theorem}[Bobkov-Ledoux \cite{BLweight}]
Assume that the following weighted logarithmic Sobolev inequality is satisfied
$$\Ent_\mu(f^2)\le 2\int|\nabla f|^2\omega d\mu.$$
Assume also that $\omega$ has a finite moment of order $p\ge 2$, then for any $\mu$-centered
1-Lipshitz function $f$, one has
$$\| f\|_p\le \sqrt{p-1}\|\omega\|_p.$$
It implies that if $\parallel \omega \parallel_p \le C$,
\begin{equation}
\mu(|f|\ge t)\le \left\{\begin{array}{ll}
2e^{-t^2/2c^2e}&{\rm if}\, 0\le t\le C\sqrt{ep}\\
2e^{-t/Ce}&{\rm if}\,C\sqrt{ep}\le t\le Cep\\
2\left(\frac{Cp}t\right)^p &{\rm if}\,Cep\le t
\end{array}\right.
\end{equation}
\end{theorem}

Remark now that the weight obtained by Bobkov-Ledoux for Cauchy
measures $\nu_\beta$ is $\omega=(\beta-1)^{-1}(1+|x|^2)^2$ whereas
ours is $\omega=C(1+|x|^2)\log(1+|x|^2)$ which thus allows
integration for $L^p(\mu)$ for a larger $p$. In addition Corollary
\ref{corcauchy} is obtained for $\beta>n/2$ instead of $\beta \geq
(n+1)/2$. Thus our result furnishes in principle a larger strip of
Gaussian concentration. However the evaluation of $C$ is quite bad
here (due mainly to the local inequality). It thus raises the
question of the optimal constant with our weight. In dimension 1,
one may use the generalized Hardy inequality.

\subsubsection{Transportation inequality.}
We give here another way to derive concentration inequality, based
on transportation inequality, as derived from logarithmic Sobolev
type inequality by Bobkov-Gentil-Ledoux \cite{BGL} using
Hamilton-Jacobi equation (see also \cite{OV} for a proof based on
PDE and optimal transport, or \cite{CatGui1} for a refined
argument). Let us give quickly the argument adapted to our setting.
First, let $d_\omega$ be the new Riemanian distance associated to
$\omega$, i.e. $C_{xy}$ is the set of all absolutely continuous
paths $\gamma:[0,1]\to\R^d$ such that $\gamma(0)=x$ and
$\gamma(1)=y$ and
$$d_\omega(x,y):=\inf_{\gamma\in C_{x,y}}\int_0^1\sqrt{\omega(\gamma(s))^{-1}\gamma'(s)^2}ds.$$
Thanks to results of Cutri-DaLio \cite{CDL} or Dragoni
\cite{dragoni} (in a more general setting, like possibly degenerate
weight), the inf-convolution $Q_t^\omega f(x):=\inf\{f(y)+\frac1t
d_\omega(x,y)\}$
 is the viscosity solution of the weighted Hamilton-Jacobi equation
 \begin{equation}
 \left\{\begin{array}{ll}
 \partial_t v +\frac12 \omega|\nabla v|^2=0&\forall(x,t)\in\R^d\times]0,\infty[,\\
 v=f&\forall (x,t)\in \R^d\times\{0\}.
 \end{array}\right.
 \end{equation}
Suppose now that $\mu$ satisfies a weighted logarithmic Sobolev  inequality with weight
$2\omega\ge1$ (the factor 2 is only for a nice formulation of the result), we apply it to the
function $f^2=e^{ tQ_t^\omega g}$ and denote $G(t)=\mu(f^2)$ so that we get, using that
$$tQ_t g=t\partial_t(tQ_t g)+\frac12|\nabla(tQ_t g)|^2$$
the differential inequality
$$tG'(t)\le G(t)\log(G(t)),\qquad \qquad G'(0)=\rho\mu(g).$$
It is now immediate to obtain that
$$\mu(e^{ Q_1g})\le e^{\mu(g)}$$
which is, by Bobkov-Goetze's result \cite{bobkov-gotze}  an equivalent formulation for a $T_2$
inequality. Summarizing this argument, we get
\begin{theorem}
Suppose that $\mu$ satisfies a weighted logarithmic Sobolev inequality with weight $2\omega$, i.e. for all nice $f$
$$\Ent_\mu(f^2)\le 2\int|\nabla f|^2\,\omega d\mu,$$
then $\mu$ satisfies the following weighted Transportation-Information inequality  ($\omega T_2$):
for all probability measure $\nu$ with $d\nu=fd\mu$
\begin{equation}\label{wTI}
W_{2,\omega}^2(\nu,\mu)\le \Ent_\mu(f).
\end{equation}
\end{theorem}
Here $W_{p,\omega}(\nu,\mu)$ is the $L^p$-Wasserstein distance
between two probability measures $\nu,\mu$ on $E$. Note that as
usual, such a ($\omega T_2$) inequality implies a $(\omega T_1)$
inequality: for all probability measure $\nu$
$$W_{1,\omega}(\nu,\mu):=\sup_{\|f\|_{Lip(\omega)}\le 1}\left(\int fd\nu-\int fd\mu\right)\le \sqrt{\Ent_\mu\left(\frac{d\nu}{d\mu}\right)}$$
where $\|f\|_{Lip(\omega)}\le 1$ means that $|f(x)-f(y)|\le d_\omega(x,y)$.

The last inequality is equivalent to the fact that for all $\mu$-centered function with
$\|f\|_{Lip(\omega)}\le 1$, $\forall r>0$,
$$\mu(|f|\ge r)\le 2 e^{-r^2/2}.$$

\subsection{Entropic convergence}

\subsubsection{The natural diffusion associated to the weighted energy.}

As is well known, logarithmic Sobolev inequality are equivalent to the exponential decay in $\L
\log \L$ of the diffusion semi group-associated to the Dirichlet form present in the inequality. We then get
that a weighted logarithmic Sobolev inequality for the measure $d\mu=e^{-V(x)}dx$
$$\Ent_\mu(f^2)\le \int |\nabla f|^2\omega d\mu$$
implies that the  semi-group $(P_t^\omega)$ with generator
$$L^\omega=\omega\Delta+(\nabla \omega-\omega\nabla V).\nabla$$
satisfies
$$\Ent_\mu(P_t^\omega f)\le e^{-t/4}\Ent_\mu(f).$$
As this semigroup is reversible with respect to $\mu$, it is certainly possible to use the results
of \cite{CGWW}, via also Lyapunov conditions, to get this convergence but it is far easier to get
a Lyapunov condition on the generator $L$ than on $L^\omega$. Note that it may also be useful when
one desires to sample from $\mu$ via a Langevin tempered diffusions type algorithm (see
\cite{DFG}): we provide here an easy way to find a diffusion coefficient leading to an exponential
entropic convergence. It has to be noted that the approach is quite different than in Hwang$\&$al \cite{HHS} or Franke$\&$al \cite{FHPS} where they add a divergence free drift to accelerate the diffusion. Moreover they are limited to cases where the initial measure $\mu$ satisfies a Poincar\'e inequality. One may also get deviation inequality for integral functional of this Markov
process, once remarked that assuming weighted logarithmic Sobolev inequality implies a
transportation cost ($\omega T_2$) inequality, then we have using once again the weighted
logarithmic Sobolev inequality: for all probability measure $\nu$ with $d\nu=fd\mu$
$$W_{2,\omega}^2(\nu,\mu)\le 2\int\frac{|\nabla f|^2}{f}\omega d\mu$$
which implies, by \cite{GLWY} that for all $\mu$-centered function $f$ with
$\|f\|_{Lip(\omega)}\le 1$ and for $(X^\omega_t)_{t\ge0}$ the Markov process with generator
$L^\Omega$: for all positive $r$
$$\P_{\nu}\left(\frac1t \int_0^tf(X^\omega_s)ds\ge r\right)\le e^{-r^2/4},$$
which may be useful in Monte-Carlo simulation.

\subsubsection{Link with weak logarithmic Sobolev inequality.}
Two of the authors with I. Gentil introduced in \cite{CGG} the weak logarithmic Sobolev
inequalities, i.e. $\mu$ satisfies (WLSI) for some non increasing function $\beta$ if for all
bounded smooth function, $\forall s>0$
\begin{equation}\label{WLSI}
\Ent_{\mu}(f^2)\le \beta(s)\int|\nabla f|^2d\mu+s Osc(f)^2 \, .
\end{equation}
This is the weak counterpart of the classical Gross logarithmic Sobolev inequalities as weak
Poincar\'e inequalities of \cite{RW} were for the usual Poincar\'e inequalities. These weak
logarithmic Sobolev inequalities are particularly useful to assert the speed of convergence
towards equilibrium (for the natural Markov process associated to $\mu$) in entropy
when dealing with particular initial measure (such as Dirac mass, not suitable to an $L^2$ analysis).\\
It was shown in \cite{CGG} that weak logarithmic Sobolev inequalities are equivalent to some
capacity/measure conditions. If in dimension one, these capacity/measure conditions can be
translated into verifiable conditions, it is no more the case in larger dimensions and only a
comparison, under some additional conditions, with Beckner inequalities (stronger than Poincar\'e)
or weak Poincar\'e inequalities  gave multidimensional examples. We will show here that weighted
logarithmic Sobolev inequalities together with some concentration estimates, enable us to obtain
weak logarithmic Sobolev inequalities, so that Lyapunov type conditions plus concentration
give a new set of conditions for weak logarithmic Sobolev inequalities.\\

\begin{theorem}
Assume that $\mu$ satisfies the following weighted logarithmic Sobolev inequality
$$\Ent_\mu(f^2)\le \int \omega |\nabla f|^2d\mu$$
then $\mu$ satisfies  a (WLSI) with function $\beta(s)=g^{-1}(s)$ where
\begin{equation}
g(r)=  \mu( B_r^c)\left[2{\mathfrak c}+\log\left(1+\frac{e^2}{\mu( B_r^c)}\right)\right]
\end{equation}
with $B_r=\{x;\,\omega\le r\}$ and ${\mathfrak c}>0$ explicit.
\end{theorem}

\begin{proof}
Let us first recall the result of Theorem 2.2 of \cite{CGG} (taking advantage of Remark 2.3), that
is a capacity measure condition for weak logarithmic Sobolev inequality.

To this end, let us recall the definition of the capacity of a given measurable set
$A\subset\Omega$:
$$Cap_\mu(A,\Omega):=\inf\left\{\int|\nabla f|^2d\mu;\,1_A\le f\le 1_\Omega\right\} $$
where the infimum is taken over all Lipschitz functions. Finally if $A$ is such that $\mu(A)<1/2$
then
$$Cap_\mu(A):=\inf\{Cap_\mu(A,\Omega);\,A\subset\Omega,\mu(\Omega)\le1/2\}.$$
A sufficient condition for (\ref{WLSI}) to hold is then: for every $A$ with $\mu(A)<1/2$,
\begin{equation}\label{capwlsi}
\forall s>0,\qquad \frac{\mu(A)\log\left(1+\frac{e^2}{\mu(A)}\right)}{\beta(s)}\le Cap_\mu(A).
\end{equation}
We cannot use directly our weighted logarithmic Sobolev inequality  with this notion of capacity
so that we introduce the natural weighted capacity
$$\overline{Cap}_\mu(A,\Omega):=\inf\left\{\int|\nabla f|^2\omega d\mu;\,1_A\le f\le 1_\Omega\right\} $$
\begin{eqnarray*}
\overline{Cap}_\mu(A)&:=&\inf\{\overline{Cap}_\mu(A,\Omega);\,A\subset\Omega,\mu(\Omega)\le1/2\}\\
&=&\inf\left\{\int|\nabla f|^2\omega d\mu;\,f:\R^d\to[0,1],f 1_A=1,\mu(f=0)\ge1/2\right\}
\end{eqnarray*}
Using Bobkov-Goetze's seminal work \cite{bobkov-gotze} or its refined version by Barthe-Roberto
\cite{bartr03sipm},  the weighted logarithmic Sobolev inequality implies that for all $A$ such
that $\mu(A)<1/2$ there exists $\mathfrak c$ such that
$$\mu(A)\log\left(1+\frac{e^2}{\mu(A)}\right)\le{\mathfrak c}\,\overline{Cap}_\mu(A).$$
Consider now the set $B_r=\{x;\,\omega\le r\}$, by a simple adaptation of the proof of Gozlan
\cite[]{gozlan2}, we get that if $A\subset B_r$
$$\overline{Cap}_\mu(A)\le 2r Cap_\mu(A)+2\mu(B_r^c).$$
Remark now that the mapping $t\to t\log(1+e^2/t)$ is concave increasing for small values  of $t$,
so that for all $A$ such that $\mu(A)\le 1/2$
\begin{eqnarray*}
\mu(A)\log\left(1+\frac{e^2}{\mu(A)}\right)&\le&\mu(A\cap B_r)\log\left(1+\frac{e^2}{\mu(A\cap B_r)}\right)+
\mu(A\cap B_r^c)\log\left(1+\frac{e^2}{\mu(A\cap B_r^c)}\right)\\
&\le& {\mathfrak c}\,\overline{Cap}_\mu(A\cap B_r)+\mu( B_r^c)\log\left(1+\frac{e^2}{\mu( B_r^c)}\right)\\
&\le& 2{\mathfrak c}r\,{Cap}_\mu(A)+\mu( B_r^c)\left[2{\mathfrak c}+\log\left(1+\frac{e^2}{\mu( B_r^c)}\right)\right].
\end{eqnarray*}
Setting $s=\mu( B_r^c)\left[2{\mathfrak c}+\log\left(1+\frac{e^2}{\mu( B_r^c)}\right)\right]$, we
conclude the proof.
\end{proof}
 If $r$ is large enough, concentration result of the previous section will give upper bounds for the second term of the left hand side.
\medskip

\subsection{Modified logarithmic Sobolev inequalities}
We will prove here that weighted logarithmic Sobolev inequalities imply modified logarithmic
Sobolev inequalities (i.e. the energy is modified). These inequalities were initially introduced
by Bobkov-Ledoux \cite{BL97}, where they show that a Poincar\'e inequality implies a logarithmic
Sobolev inequality for a particular class of functions ($|\nabla f/f|\le c<C_{c_{SG}}$ where
$c_{SG}$ is the spectral gap constant). These results were later extended to measures between
exponential and Gaussian by Gentil and al \cite{GGM1,GGM2}. For recent results, giving nice
conditions we will discuss later, see also \cite{BK08}.

\begin{theorem}
Let $H$ and $H^y$ be a pair of dual convex  Young functions, such that $H(|x|)/|x|\ge a>0$
for large $|x|$ and $H^*(\epsilon |x|)\le b(\epsilon) H^y(|x|)$ with $b(\epsilon)\to0$ as $\epsilon\to0$.\\
Suppose now that the following weighted logarithmic Sobolev inequality holds
$$\Ent_\mu(f^2)\le \int |\nabla f|^2 \,\omega d\mu$$
for some weight $\omega\ge1$, that a Poincar\'e inequality holds and that for some $\alpha>0$
\begin{equation}\label{integcond}
K:=\int e^{\alpha H^y(\omega)}d\mu<\infty.
\end{equation}
Then the following modified logarithmic Sobolev inequality holds
\begin{equation}\label{mlsi}
\Ent_\mu(f^2)\le C\int\left(H\left(\epsilon^{-1}\left|\frac{\nabla f}f\right|^2\right)f^2+|\nabla
f|^2\right)d\mu
\end{equation}
for sufficiently small $\epsilon$ and some constant $C$ (explicit in the proof).
\end{theorem}

\begin{proof}
Actually, it is sufficient to get a defective modified logarithmic Sobolev inequality, since a
Poincar\'e inequality allows us  to tighten a defective inequality thanks to \cite[Th. 2.4]{BK08}.
We then have
\begin{eqnarray*}
\Ent_\mu(f^2)&\le& \int |\nabla f|^2 \,\omega d\mu\\
&=& \int \epsilon^{-1}\left|\frac{\nabla f}f\right|^2 \,\epsilon\omega f^2 d\mu\\
&\le&\int H\left(\epsilon^{-1}\left|\frac{\nabla f}f\right|^2\right)f^2d\mu+\int
H^y(\epsilon\omega)f^2d\mu.
\end{eqnarray*}
Choose now $\epsilon$ sufficiently small so that $b(\epsilon)\le \alpha/2$ so that
\begin{eqnarray*}
\int H^*(\epsilon\omega)f^2d\mu&\le&\frac12\int\alpha H^y(\omega)f^2d\mu\\
&\le& \frac12\int (\alpha H^y(\omega)-\log K)f^2d\mu+\frac12\log K\int f^2d\mu\\
&\le&\frac12\Ent_\mu(f^2)+\frac12\log K\int f^2d\mu
\end{eqnarray*}
where we have used the variational formula for the entropy in the last line. Plugging the latter
inequality in the preceding one, we obtain the defective modified logarithmic Sobolev inequality:
$$\Ent_\mu(f^2)\le 2\int H\left(\epsilon^{-1}\left|\frac{\nabla f}f\right|^2\right)f^2d\mu+\log(K)\int f^2d\mu$$
which ends the proof.
\end{proof}
One may then use the Lyapunov conditions used to derive a weighted logarithmic Sobolev inequality
to get a generalization of Barthe-Kolesnikov \cite[Th. 5.27,5.28]{BK08}.\\ \medskip

 {\bf Examples: }\\
Consider the usual (for modified LSI) examples: $d\mu=Z_\alpha e^{-|x|^\beta}$ for $1<\alpha\le2$
so that the Poincar\'e inequality is valid. Using Lyapunov function $W(x)=e^{a|x|^\beta}$ for $a$
less than one, one may easily derive the following Lyapunov condition:
$$LW\le -c|x|^{2(\beta-1)}W+b1_{B(0,R)}$$ from which one deduces using $\psi(w)=\log(w)$ and Theorem \ref{mainTH} (and Prop. \ref{eqFsob}):
$$\Ent_\mu(f^2)\le C\int |\nabla f|^2 \,(1+|x|^{2-\beta})d\mu.$$
Consider now the Young functions $H_\beta(x)=|x|^{\beta\over 2(\beta-1)}$ and
$H_\beta^y(x)=c_\beta|x|^{\beta\over 2-\beta}$ so that $H_\beta^y(\epsilon\omega)=c_\beta
\omega^{\beta\over 2-\beta}|x|^\beta$ which is easily seen to be integrable wrt $\mu$ for
$\epsilon$ sufficiently small. We then get
 $$\Ent_\mu(f^2)\le C\int\left(\left|\frac{\nabla f}f\right|^{\beta\over\beta-1}f^2+|\nabla f|^2\right)d\mu$$
 for some constant $C$, which is a generalization in the multidimensional case of \cite{GGM1}.
\bigskip

\bibliographystyle{plain}

\begin{thebibliography}{10}

\bibitem{bakry}
D.~Bakry.
\newblock L'hypercontractivit\'e et son utilisation en th\'eorie des
  semigroupes.
\newblock In {\em Lectures on Probability theory. \'Ecole d'{\'e}t{\'e} de
  {P}robabilit{\'e}s de St-Flour 1992}, volume 1581 of {\em Lecture Notes in
  Math.}, pages 1--114. Springer, Berlin, 1994.

\bibitem{BBCG}
D.~Bakry, F.~Barthe, P.~Cattiaux, and A.~Guillin.
\newblock A simple proof of the {P}oincar\'e inequality for a large class of
  probability measures.
\newblock {\em Electronic Communications in Probability.}, 13:60--66, 2008.

\bibitem{BCG}
D.~Bakry, P.~Cattiaux, and A.~Guillin.
\newblock Rate of convergence for ergodic continuous {M}arkov processes :
  {L}yapunov versus {P}oincar\'e.
\newblock {\em J. Func. Anal.}, 254:727--759, 2008.

\bibitem{BCR1}
F.~Barthe, P.~Cattiaux, and C.~Roberto.
\newblock Interpolated inequalities between exponential and {G}aussian,
  {O}rlicz hypercontractivity and isoperimetry.
\newblock {\em {R}ev. {M}at. {I}ber.}, 22(3):993--1066, 2006.

\bibitem{BK08}
F.~Barthe and A.~V. Kolesnikov.
\newblock Mass transport and variants of the logarithmic {S}obolev inequality.
\newblock {\em J. Geom. Anal.}, 18(4):921--979, 2008.

\bibitem{bartr03sipm}
F.~Barthe and C.~Roberto.
\newblock Sobolev inequalities for probability measures on the real line.
\newblock {\em Studia Math.}, 159(3), 2003.

\bibitem{BZ}
F.~Barthe and Z.. Zhang.
\newblock Private communication.
\newblock 2009.

\bibitem{blanchet}
A.~Blanchet, M.~Bonforte, J.~Dolbeault, G.~Grillo, and J.L. V{\'a}zquez.
\newblock Asymptotics of the fast diffusion equation via entropy estimates.
\newblock {\em Arch. Ration. Mech. Anal.}, 191(2):347--385, 2009.

\bibitem{BGL}
S.~G. Bobkov, I.~Gentil, and M.~Ledoux.
\newblock Hypercontractivity of {H}amilton-{J}acobi equations.
\newblock {\em J. Math. Pu. Appli.}, 80(7):669--696, 2001.

\bibitem{bobkov-gotze}
S.~G. Bobkov and F.~G{\"o}tze.
\newblock Exponential integrability and transportation cost related to
  logarithmic {S}obolev inequalities.
\newblock {\em J. Funct. Anal.}, 163(1):1--28, 1999.

\bibitem{BL97}
S.G. Bobkov and M.~Ledoux.
\newblock Poincar\'e inequalities and {T}alagrand concentration phenomenon for
  the exponential distribution.
\newblock {\em Prob. Theor. Rel. Fields}, 107:383--400, 1997.

\bibitem{BLweight}
S.G. Bobkov and M.~Ledoux.
\newblock Weighted {P}oincar\'e-type inequalities for {C}auchy and other convex
  measures.
\newblock {\em Ann. Probab.}, 37(2):403--427, 2009.

\bibitem{cat4}
P.~Cattiaux.
\newblock A pathwise approach of some classical inequalities.
\newblock {\em Potential Analysis}, 20:361--394, 2004.

\bibitem{CGG}
P.~Cattiaux, I.~Gentil, and A.~Guillin.
\newblock Weak logarithmic {S}obolev inequalities and entropic convergence.
\newblock {\em Probab. Theory Related Fields}, 139(3-4):563--603, 2007.

\bibitem{CGGR}
P.~Cattiaux, N.~Gozlan, A.~Guillin, and C.~Roberto.
\newblock Functional inequalities for heavy tailed distributions and
  applications to isoperimetry.
\newblock available on {M}ath. {A}r{X}iv. 0807.3112. [math {PR}], 2008.

\bibitem{CatGui1}
P.~Cattiaux and A.~Guillin.
\newblock On quadratic transportation cost inequalities.
\newblock {\em J. Math. Pures Appl.}, 88(4):341--361, 2006.

\bibitem{CGgre}
P.~Cattiaux and A.~Guillin.
\newblock Functional inequalities via {L}yapunov conditions.
\newblock to appear in Proceedings of the summer school on Optimal Transport
  (Grenoble 2009). available on {M}ath. {A}r{X}iv. 1001.1822. [math {PR}],
  2010.

\bibitem{CGWW}
P.~Cattiaux, A.~Guillin, F-Y. Wang, and L.~Wu.
\newblock Lyapunov conditions for super {P}oincar\'e inequalities.
\newblock {\em J. Funct. Anal.}, 256(6):1821--1841, 2009.

\bibitem{CGW}
P.~Cattiaux, A.~Guillin, and L.~Wu.
\newblock A note on {T}alagrand transportation inequality and logarithmic
  {S}obolev inequality.
\newblock To appear in {Prob. Theo. Rel. Fields}., 2008.

\bibitem{CDL}
A.~Cutr{\`{\i}} and F.~Da~Lio.
\newblock Comparison and existence results for evolutive non-coercive
  first-order {H}amilton-{J}acobi equations.
\newblock {\em ESAIM Control Optim. Calc. Var.}, 13(3):484--502, 2007.

\bibitem{DMC}
J.~Denzler and R.~J. McCann.
\newblock Fast diffusion to self-similarity: complete spectrum, long-time
  asymptotics and numerology.
\newblock {\em Arch. Ration. Mech. Anal.}, 175(3):301--342, 2005.

\bibitem{DV3}
M.~D. Donsker and S.~R.~S. Varadhan.
\newblock Asymptotic evaluation of certain {M}arkov process expectations for
  large time. {III}.
\newblock {\em Comm. Pure Appl. Math.}, 29(4):389--461, 1976.

\bibitem{DV4}
M.~D. Donsker and S.~R.~S. Varadhan.
\newblock Asymptotic evaluation of certain {M}arkov process expectations for
  large time. {IV}.
\newblock {\em Comm. Pure Appl. Math.}, 36(2):183--212, 1983.

\bibitem{DFG}
R.~Douc, G.~Fort, and A.~Guillin.
\newblock Subgeometric rates of convergence of {$f$}-ergodic strong {M}arkov
  processes.
\newblock {\em Stochastic Process. Appl.}, 119(3):897--923, 2009.

\bibitem{dragoni}
F.~Dragoni.
\newblock Metric {H}opf-{L}ax formula with semicontinuous data.
\newblock {\em Discrete Contin. Dyn. Syst.}, 17(4):713--729, 2007.

\bibitem{FHPS}
B.~Franke, C.R.~Hwang, H.M.~Pai, and S.J.~Sheu.
\newblock The behavior of the spectral gap under growing drift.
\newblock{\em Trans. Amer. Math. Soc}, 362(3):1325--1350, 2009.

\bibitem{GGM1}
I.~Gentil, A.~Guillin, and L.~Miclo.
\newblock Modified logarithmic {S}obolev inequalities and transportation
  inequalities.
\newblock {\em Probab. Theory Related Fields}, 133(3):409--436, 2005.

\bibitem{GGM2}
I.~Gentil, A.~Guillin, and L.~Miclo.
\newblock Modified logarithmic sobolev inequalities in null curvature.
\newblock To appear in {R}evista {M}atematica {I}beroamericana, 2006.

\bibitem{gozlan2}
N.~Gozlan.
\newblock Poincar\'e inequalities for non euclidean metrics and transportation
  inequalities.
\newblock Preprint. {Available on Math ArXiv 0707.2834[ps]}, 2007.

\bibitem{G2}
A.~Guillin.
\newblock Averaging principle of {SDE} with small diffusion: moderate
  deviations.
\newblock {\em Ann. Probab.}, 31(1):413--443, 2003.

\bibitem{G1}
Arnaud Guillin.
\newblock Moderate deviations of inhomogeneous functionals of {M}arkov
  processes and application to averaging.
\newblock {\em Stochastic Process. Appl.}, 92(2):287--313, 2001.

\bibitem{GLWY}
A.~Guillin, C.~L{\'e}onard, L.~Wu, and N.~Yao.
\newblock Transportation-information inequalities for {M}arkov processes.
\newblock {\em Probab. Theory Related Fields}, 144(3-4):669--695, 2009.

\bibitem{HHS}
C.R.~Hwang, S.Y.~Hwang-Ma, and S.J~Sheu.
\newblock Accelerating diffusions.
\newblock {\em Ann. Appl. Prob.}, 15:1433--1444, 2005/

\bibitem{KM1}
I.~Kontoyiannis and S.~P. Meyn.
\newblock Spectral theory and limit theorems for geometrically ergodic {M}arkov
  processes.
\newblock {\em Ann. Appl. Probab.}, 13(1):304--362, 2003.

\bibitem{KM2}
I.~Kontoyiannis and S.~P. Meyn.
\newblock Large deviations asymptotics and the spectral theory of
  multiplicatively regular {M}arkov processes.
\newblock {\em Electron. J. Probab.}, 10:no. 3, 61--123 (electronic), 2005.

\bibitem{MT}
S.~P. Meyn and R.~L. Tweedie.
\newblock {\em Markov chains and stochastic stability}.
\newblock Communications and Control Engineering Series. Springer-Verlag London
  Ltd., London, 1993.

\bibitem{MT2}
S.~P. Meyn and R.~L. Tweedie.
\newblock Stability of markovian processes {II}: continuous-time processes and
  sampled chains.
\newblock {\em Adv. Appl. Proba.}, 25:487--517, 1993.

\bibitem{MT3}
S.~P. Meyn and R.~L. Tweedie.
\newblock Stability of markovian processes {III}: {F}oster-{L}yapunov criteria
  for continuous-time processes.
\newblock {\em Adv. Appl. Proba.}, 25:518--548, 1993.

\bibitem{OV}
F.~Otto and C.~Villani.
\newblock Generalization of an inequality by {T}alagrand and links with the
  logarithmic {S}obolev inequality.
\newblock {\em J. Funct. Anal.}, 173(2):361--400, 2000.

\bibitem{RW}
M.~R{\"o}ckner and F.~Y. Wang.
\newblock Weak {P}oincar\'e inequalities and {$L\sp 2$}-convergence rates of
  {M}arkov semigroups.
\newblock {\em J. Funct. Anal.}, 185(2):564--603, 2001.

\bibitem{w00}
F.~Y. Wang.
\newblock Functional inequalities for empty essential spectrum.
\newblock {\em J. Funct. Anal.}, 170(1):219--245, 2000.

\bibitem{Wbook}
F.~Y. Wang.
\newblock {\em Functional inequalities, {M}arkov processes and {S}pectral
  theory}.
\newblock Science Press, Beijing, 2005.

\bibitem{wu1}
Liming Wu.
\newblock Large and moderate deviations and exponential convergence for
  stochastic damping {H}amiltonian systems.
\newblock {\em Stochastic Process. Appl.}, 91(2):205--238, 2001.

\end{thebibliography}

\def\cprime{$'$}

\end{document}